\documentclass[10pt]{amsart}
\usepackage{amssymb,amsmath,amscd,amsthm,pb-diagram}
\usepackage[matrix,arrow]{xy}

\title{Manifolds homotopy equivalent to $P^n \# P^n$ }
\author{Jeremy Brookman}
\author{James F. Davis}
 \author{Qayum Khan}
\thanks{Partially supported by a grant from the National Science Foundation}
\address{Department of Mathematics \\
Indiana University \\
Bloomington, IN 47405 \\ USA} \email{jeremy@jbrookman.me.uk,
jfdavis@indiana.edu, qkhan@indiana.edu}
\date{}

\newcommand{\cS}{{\mathcal{S}}}
\newcommand{\cL}{{\mathcal{L}}}
\newcommand{\cN}{{\mathcal{N}}}

\newcommand{\RP}{{P}}
\newcommand{\C}{{\mathbb C}}
\newcommand{\F}{{\mathbb F}}
\newcommand{\Q}{{\mathbb Q}}
\newcommand{\R}{{\mathbb R}}
\newcommand{\Z}{{\mathbb Z}}
\newcommand{\PL}{{\mathrm{PL}}}
\newcommand{\even}{\mathrm{ev}}
\newcommand{\odd}{\mathrm{od}}
\newcommand{\inv}{^{-1}}
\newcommand{\iso}{\cong}
\newcommand{\xra}{\xrightarrow}
\newcommand{\Inn}[2]{\langle #1, #2 \rangle}
\newcommand{\prn}[1]{\left( #1 \right)}
\newcommand{\brk}[1]{\left[ #1 \right]}
\newcommand{\SmMatrix}[1]{\left(\begin{smallmatrix}#1\end{smallmatrix}\right)}
\newcommand{\SmBMatrix}[1]{\left[\begin{smallmatrix}#1\end{smallmatrix}\right]}
\newcommand{\Span}[1]{\mathrm{span}\left\{ #1 \right\}}
\newcommand{\ol}[1]{\bar{#1}}

\newcommand{\wt}[1]{\widetilde{#1}}
\DeclareMathOperator{\aqk}{{aqk}}
 
\DeclareMathOperator{\haut}{{hAut}}

\DeclareMathOperator{\out}{{Out}}
 \DeclareMathOperator{\id}{{id}}
\DeclareMathOperator{\nncc}{{nncc}}  \DeclareMathOperator{\Split}{{split}}
 \DeclareMathOperator{\sw}{{sw}}
\DeclareMathOperator{\unil}{{UNil}} \DeclareMathOperator{\Wh}{{Wh}}

\newtheorem{theorem}{Theorem}

\newtheorem{corollary}[theorem]{Corollary}
\newtheorem{lemma}[theorem]{Lemma}
\theoremstyle{definition}
\newtheorem{definition}[theorem]{Definition}
\newtheorem{remark}[theorem]{Remark}
\newtheorem{example}[theorem]{Example}

\begin{document}

\begin{abstract}  We classify, up to homeomorphism, all closed manifolds having the homotopy type of a connected sum of two copies of real projective $n$-space.
\end{abstract}

\maketitle

\section{Statement of results}

Let $P^n = P_n(\R)$ be real projective $n$-space. L\'opez de Medrano \cite{LdM} and C.T.C. Wall \cite{WallBAMS,Wall} classified, up to $\PL$
homeomorphism, all closed $\PL$ manifolds homotopy equivalent to $\RP^n$ when $n > 4$. This was extended to the topological category by
Kirby-Siebenmann \cite[p. 331]{KS}. Four-dimensional surgery \cite{FQ} extends the homeomorphism classification to dimension 4.

Cappell \cite{CappellFree, CappellUnitary,CappellSplitting} discovered that the situation for connected sums is much more complicated.  In
particular, he showed \cite{CappellConnected} that there are closed manifolds homotopy equivalent to $\RP^{4k+1} \# \RP^{4k+1}$ which are not
non-trivial connected sums.  Recent computations of the unitary nilpotent group for the integers by Connolly-Ranicki \cite{CR}, Connolly-Davis
\cite{CD}, and Banagl-Ranicki \cite{BR} show that there are similar examples in dimension $4k$ (see \cite{JK} for an analysis when $k=1$).

In this paper we classify up to homeomorphism all closed manifolds homotopy equivalent to $\RP^n \# \RP^n$.  Any such manifold has $S^{n-1}
\times S^1$ as a two-fold cover; equivalently we classify free involutions on $S^{n-1} \times S^1$ inducing a non-trivial map on $H_1$.

This paper was prompted by a question of Wolfgang L\"uck \cite[Sequence (4.10), Theorem 4.11]{Lueck} -- what does the group automorphism $\Z_2 *
\Z_2 \to \Z_2
* \Z_2$ given by interchanging the $\Z_2$'s induce in $L$-theory?  We give a complete answer to L\"uck's question and apply the answer to the
classify the above manifolds.

By the positive solution to Kneser's conjecture (see \cite{Hempel}) any closed 3-manifold homotopy equivalent to $P^3 \# P^3$ is homeomorphic to $Q^3 \# R^3$ where $Q^3$ and $R^3$ are closed 3-manifolds homotopy equivalent to $P^3$.  The spherical space form conjecture in dimension 3 would imply that $Q^3$ and $R^3$ are homeomorphic to $P^3$.  Hence, conjecturally,  any closed 3-manifold homotopy equivalent to $P^3 \# P^3$ is standard.  Henceforth we assume $n > 3$.  Note that the fundamental group of $P^4 \# P^4$ is small in the sense of Freedman-Quinn \cite{FQ}, so that surgery theory applies.

Let $\bar I_n$ (respectively $\bar J_n$) be the set of homeomorphism classes of closed manifolds homotopy equivalent to $\RP^n$ (respectively
$\RP^n \# \RP^n$).  For $n$ even, let $I_n = \bar I_n$ and $J_n = \bar J_n$.  For $n$ odd, let $I_n$ (respectively $J_n$) be the set of oriented
homeomorphism classes of closed oriented manifolds homotopy equivalent to $\RP^n$ (respectively $\RP^n \# \RP^n$).  The set $I_n$ was computed
in \cite[\S 16 Annex 3]{KS}; we review the computation in Section \ref{structure_set}.

Let $R$ be a ring with involution.  Let $A$ and $B$ be $(R,R)$-bimodules with involution.  (The case of interest is $R =\Z$, and then $A$ and
$B$ are just abelian groups with an automorphism of order 2.)  Cappell \cite{CappellUnitary} defined unitary nilpotent groups $\unil_n(R;A,B)$,
abelian groups with the following properties.
\begin{itemize}

\item  They satisfy periodicity $\unil_n(R; A,B) \cong \unil_{n+4}(R;A,B)$ and semiperiodicity $\unil_n(R; A,B) \cong \unil_{n+2}(R;A^-,B^-)$, where $A^-$ is the original bimodule but with the involution $a \mapsto \bar a$ replaced by $a \mapsto -\bar a$.

\item  They obstruct splitting:  Given a  homotopy equivalence $h : X \to \RP^n \# \RP^n$ for any $n>4$, there is an element $\Split(h) \in \unil_{n+1}(\Z;\Z^\epsilon,\Z^\epsilon)$ which vanishes if and only if $h$ is splittable\footnote{$h$ is splittable if it is homotopic to a map $f$, transverse to $S$, which induces a homotopy equivalence $f^{-1}(S) \to S$,  where $S$ is the codimension one sphere defining the connected sum.}.  Here $\epsilon = (-1)^{n+1}$.  Conversely all obstructions are realized.

Note that to define the element $\Split(h)$ one needs to choose an orientation for $P^n \# P^n$, that is, a generator for the infinite cyclic group $H_n(P^n \# P^n; \Z^w)$ where $w$ is the orientation character.

\item  There is a split injection $i: \unil_n(\Z;\Z^\epsilon,\Z^\epsilon) \to L_{n}(\Z[\Z_2^\epsilon * \Z_2^\epsilon])$.

\end{itemize}

The groups $\unil_n(\Z;\Z,\Z)$ have been completely computed (see, for example,  \cite{CD}).   The groups $\unil_0(\Z;\Z,\Z)$ and $\unil_1(\Z;\Z,\Z)$ vanish.  As abelian groups 
$\unil_2(\Z;\Z,\Z)  \cong \oplus_{\infty} \Z_2$ and $\unil_3(\Z;\Z,\Z)  \cong \oplus_{\infty} \Z_2 \oplus_{\infty} \Z_4$.

The ``switch map'' $\sw: \RP^n \# \RP^n \to \RP^n \# \RP^n$ interchanges the two summands.  We also use the same notation to denote the induced
involution on $\Z_2 * \Z_2$, on $L_{n}(\Z[\Z_2^\epsilon * \Z_2^\epsilon])$, and on $\unil_n(\Z;\Z^\epsilon,\Z^\epsilon)$.  The semiperiodicity mentioned above is equivariant with respect to the switch map.

The following theorem, proved in Sections  \ref{unil_2} and  \ref{unil_3} is the main technical result of this paper.

\begin{theorem} \label{sw}
\begin{enumerate}
\item  The map $\sw : \unil_2(\Z; \Z,\Z) \to \unil_2(\Z; \Z,\Z)$ is the identity.
\item  Connolly-Davis \cite{CD} give a $\Z$-module isomorphism
$$\SmMatrix{j_1 & j_2}:\frac{t\Z_4[t]}{\{2p(t^2)-2p(t) : p(t) \in t\Z_4[t] \}}\times t\Z_2[t]\to \unil_3(\Z;\Z,\Z).$$
Let $\pi:\frac{t\Z_4[t]}{\{2p(t^2)-2p(t)\}}\to t\Z_2[t]$ be the quotient map $\pi[tp] = [tp]$.
Then using the above coordinates $\SmMatrix{j_1 & j_2}$ for $\unil_3(\Z;\Z,\Z)$, the switch map is $\sw = \SmMatrix{1 & 0\\ \pi & 1}$.
\end{enumerate}

\end{theorem}

Thus the  switch map $\sw : \unil_3(\Z; \Z, \Z) \to \unil_3(\Z; \Z, \Z)$
  is the identity 
on all elements which are multiples of two, as well as on some elements which are not multiples of two, but is non-trivial on some elements of order two.

Given homotopy equivalences $h_1: X_1 \to \RP^n$ and $h_2: X_2 \to \RP^n$ and
an element $\vartheta \in \unil_{n+1}(\Z;\Z^\epsilon,\Z^\epsilon)$, Wall realization \cite[Theorem 10.4 and 10.5]{Wall} produces a normal bordism 
$$
(g; h_1 \# h_2, h) : (W; X_1 \# X_2,X) \to P^n \# P^n \times ([0,1]; \{0\}, \{1\}) 
$$
with $h : X \to P^n \# P^n$ a homotopy equivalence and with the rel $\partial$ surgery obstruction of $g$ the image of $\vartheta$ in the $L$-group
$$
\sigma_*(g) = i(\vartheta) \in i(\unil_{n+1}(\Z;\Z^\epsilon,\Z^\epsilon)) \subset L_{n+1}(\Z[\Z_2^\epsilon * \Z_2^\epsilon])
$$

 Section \ref{structure_set} proves the
following theorem.

\begin{theorem} \label{main}
Let $n > 3$ and $\epsilon = (-1)^{n+1}$. There is a bijection from the set
\[
\prn{\text{unordered pairs from the set } I_n} \times \prn{\unil_{n+1}(\Z;\Z^\epsilon,\Z^\epsilon) / \sw} \to J_n
\]
defined by
\[
\prn{\{[h_1: X_1 \to \RP^n], [h_2: X_2 \to \RP^n]\}, [\vartheta]}
\mapsto [X],
 \]
 where the homotopy equivalence $h : X \to \RP^n\# \RP^n$
 is produced by Wall realization as above.

Moreover if $[\vartheta] \neq 0$, then $[X]$ is not represented by
the connected sum of manifolds with fundamental group $\Z_2$.

For $n \not \equiv 3 \pmod 4$, $I_n \to \bar I_n$ and $J_n \to \bar J_n$ are bijective, while for $n \equiv 3 \pmod 4$, $I_n \to \bar I_n$ and
$J_n \to \bar J_n$ are at most 2-to-1.

\end{theorem}

Theorem \ref{sw} computes $\unil_n(\Z;\Z,\Z)/\sw$ for $n=2,3$; the low-dimensional manifolds to keep in
mind are $\RP^5 \# \RP^5$ and $\RP^4 \# \RP^4$, respectively.  A precise description of the maps $I_n \to \bar I_n$ and $J_n \to \bar J_n$ is given at the end of Section \ref{structure_set}.

The computation of the switch map on UNil should be considered as the main result of this paper.  We now indicate the difficulty.  For a ring
with involution $R$, define $NL_n(R)$ to be the kernel of the augmentation map $L_n(R[t]) \to L_n(R)$ given by sending $t$ to $0$.  The papers \cite{CD}, \cite{CR}, \cite{BR} proceed to compute $NL_n(R)$, and hence $\unil_n(\Z;\Z,\Z)$ since
Connolly-Ranicki \cite[Definition 13, Theorem A]{CR} define a map
$$
r : \unil_n(R;R,R) \to NL_n(R)
$$
and show it is an isomorphism.  However the induced switch map on $NL_n(R)$ is not apparent and the inverse map $r^{-1}$ is not explicit.

Here is how to get around the difficulty.  Cappell \cite{CappellUnitary} defined a split injection $i:\unil_n(\Z;\Z,\Z)\to L_n(\Z[\Z_2 *\Z_2])$ which commutes with the switch map induced by the ring automorphism of $\Z[\Z_2 *\Z_2]$ given by switching factors.  We give an explicit formula for $F=i\circ r^{-1}$.  Then for $x$ in a certain generating set of $NL_n(\Z)$, we find $y$ such that $\sw F(x)=F(y)$, and thereby deduce that $\sw(r^{-1}(x)) = r^{-1}(y)$.

\section{The structure set of $\RP^n \# \RP^n$} \label{structure_set}

In this section we compute the structure set of $\RP^n \# \RP^n$ and
reduce the homeomorphism classification to the structure set of
$\RP^n \# \RP^n$ modulo the $\Z_2$-action given by the switch map.
The material in the section is a standard application of the surgery
machine.

The structure set $\cS(M)$ of a closed topological manifold $M$ is
defined to be the set of equivalence classes of ``$s$-triangulations,''
simple homotopy equivalences $h : X \to M$ where $X$ is a closed
topological manifold. Two such $s$-triangulations $h_1 : X_1 \to M$
and $h_2 : X_2 \to M$ are equivalent if there is a homeomorphism $f
: X_1 \to X_2$ so that $h_2 \circ f$ and $h_1$ are homotopic.

For a space $M$ having the homotopy type of a finite CW complex, let $\haut(M)$ be the group of homotopy classes of simple self-homotopy equivalences $M \to M$.  The following lemma is
standard and its proof is trivial.

\begin{lemma}\label{Lemma_Reduction}
$\haut(M)$ acts on $\cS(M)$ by post-composition.  The forgetful map
\begin{align*}
\cS(M)/\haut(M) & \to \frac{\text{closed manifolds simple homotopy equivalent to $M$}}{\text{homeomorphism}}\\
[h : X \to M] & \mapsto [X]
\end{align*}
is a bijection of sets.
\end{lemma}

Recall the fundamental group of $\RP^n\#\RP^n$ for $n > 2$ is isomorphic to both  the infinite dihedral group,
\[
D_\infty =\Z \rtimes \Z_2 = \langle t, a \;|\; ata = t\inv, a^2 = 1\rangle \]
 and to the free product (letting $t=ba$)
 \[
 \Z_2 * \Z_2 = \langle a, b \;|\; a^2 = 1 =
b^2\rangle, \] which has \cite[Theorem 6.4]{Stallings} Whitehead torsion group $\Wh(\Z_2) \oplus \Wh(\Z_2) = 0$.  Thus for
$M = P^n \# P^n$ or $P^n$ we can drop the word ``simple'' in the definition of $\cS(M)$ and $\haut(M)$.

We next review the computations of $\cS(P^n)$, $I_n$, and $\bar I_n$.  Write $n = 4m + \ell > 3$ for unique $m \geq 0$ and $0 < \ell \leq 4$. Recall that \cite[p. 331]{KS} the topological structure set
$\cS(\RP^n)$ is in bijection with the set
\[
\bigoplus_{0 < k \leq 2m + [\ell/4]} \Z_2 \oplus \begin{cases}\Z & \text{if } \ell=3\\ 0 & \text{if } \ell\not =3.\end{cases}
\]
The $\Z_2$ summands (normal invariants) arise as the surgery obstructions of degree one normal maps which are the transverse restrictions to $\RP^{2k}$ of the given homotopy equivalence to $\RP^n$.  The map from $\Z$ to the structure set is given by Wall realization $\Z \cong \widetilde L_{n+1}(\Z[\Z_2]) \to \cS(P^n)$ and the splitting map from the structure set to $\Z$ is the Browder-Livesay desuspension invariant of a free $\Z_2$-action on $S^{4m+3}$ to an action on some embedded
$S^{4m+2}$.

It is (almost) elementary to show that any self-homotopy equivalence of $\RP^n$ which does not reverse orientation is homotopic to the identity.
Thus for $n$ even, $\haut(\RP^n)= *$ and for $n$ odd, $\haut(\RP^n) = \Z_2$.  This acts trivially on the normal invariant, but for $n = 4m +3$
acts by multiplication by $-1$ on the $\Z$-summand, since reversing orientation reverses the sign of any signature invariant.  The following
theorem summarizes the discussion.

\begin{theorem} \label{P_classify} Let  $n = 4m + \ell > 3$ where $m \geq 0$ and $0 < \ell \leq 4$.
There are bijections
\begin{align*}
\Phi : \cS(P^n) &\longrightarrow \bigoplus_{0 < k \leq 2m + [\ell/4]} \Z_2 \oplus \begin{cases}\Z & \text{if } \ell=3\\ 0 & \text{else}.\end{cases}\\
I_n &\longrightarrow \cS(P^n) \\
\bar I_n & \longrightarrow \frac{ \cS(P^n)}{\Phi^{-1}(z) \sim \Phi^{-1}(-z)}
\end{align*}
\end{theorem}

Identify
$$
\RP^n \# \RP^n = \frac{S^{n-1} \times S^1}{ (w,z) \sim (-w,\bar z)}
$$
The sphere defining the connected sum is $S^{n-1} \times \{\pm i \}/\sim$.
 Following Cappell \cite[Proof 3]{CappellConnected}, define
self-homeomorphisms $\gamma_1,\gamma_2,\gamma_3$ of $\RP^n \# \RP^n$
by:

\begin{itemize}
\item $\gamma_1[w,z] = [w,-z]$,
which interchanges the two summands of $\RP^n \# \RP^n$; also called the switch map $\sw$,
\item $\gamma_2[(w_1, w_2, \dots, w_n),z] = [(w_1, w_2, \dots, w_{n-1}, -w_n),z]$,
which reflects through $\RP^{n-1} \# \RP^{n-1}$, and
\item $\gamma_3[w,z] = \begin{cases}
[\tau(z^2)(w),z] & \text{ if } \mathrm{Im}\,z \geq 0\\
[\tau(\bar z^2)(w),z] & \text{ if } \mathrm{Im}\,z \leq 0
\end{cases}$, which Dehn twists along the connecting cylinder via the isotopy $\tau: S^1 \to SO(n)$ generating $\pi_1(SO(n)) \iso \Z_2$.
\end{itemize}

When $n$ is even, $\gamma_2$ is isotopic to the identity via $$[(w_1 \cos (\pi t) - w_2 \sin(\pi t), w_1 \sin (\pi t) + w_2 \cos (\pi t), \dots , -w_n),z)].$$

\begin{lemma}[{\cite[Proof 3]{CappellConnected}}, see also \cite{JK}]\label{Lemma_hAut} \label{gen}
  $\haut(\RP^n \# \RP^n) = \langle \gamma_1, \gamma_2, \gamma_3 \rangle$.  For $n$ even, $\haut(\RP^n \# \RP^n) = \langle \gamma_1, \gamma_3 \rangle$.
\end{lemma}

In particular, all elements of $\haut(\RP^n \# \RP^n)$ are splittable.

\begin{corollary}[{\cite[Lemma 2]{CappellConnected}}]
 \label{Cor_NotSplittable}
If an $s$-triangulation $h : X \to \RP^n \# \RP^n$ is not
splittable, then $X$ is not a connected sum of manifolds with
fundamental group $\Z_2$.
\end{corollary}

\begin{lemma}\label{Lemma_SplitAction}
Let $\cS_{\Split}(P^n \# P^n)$ be the subset of $\cS(P^n \# P^n)$ given by splittable homotopy equivalences.  Let $n > 3$ and $\epsilon = (-1)^{n+1}$.
\begin{enumerate}
\item  Cappell's nilpotent normal cobordism construction gives a bijection
\begin{align*}
\nncc : \cS(P^n \# P^n) &\to \cS_{\Split}(P^n \# P^n) \times \unil_{n+1}(\Z;\Z^\epsilon,\Z^\epsilon)\\
h & \mapsto (-\Split(h)\cdot h,\Split(h)),
\end{align*}
 where $\cdot $ refers to the action of the $L$-group on the structure set.

\item The above bijection is equivariant with to the action of $\haut(P^n \# P^n)$ on the
structure and split structure sets given by post-composition and  on $ \unil_{n+1}(\Z;\Z^\epsilon,\Z^\epsilon)$ given by the switch map if the homotopy automorphism represents the non-trivial element of $\out(\Z_2 * \Z_2)$ and by the identity otherwise.

\item  Connected sum gives a bijection
$\# : \cS(P^n) \times \cS(P^n) \to \cS_{\Split}(P^n \# P^n)$.

\end{enumerate}

\end{lemma}

\begin{proof}
(1) The action of the $L$-group on the structure set is given by Wall realization, see \cite[Theorem 10.4 and 10.5]{Wall}.

 For $n > 5$, the nilpotent normal cobordism construction is given in \cite[\S II.1, \S III.2]{CappellSplitting} (see also  \cite[Theorem 3]{CappellFree}).  In \cite[Theorem 1]{CappellFree}  the map $\Split : \cS(P^n \# P^n) \to
 \unil_{n+1}(\Z;\Z^\epsilon,\Z^\epsilon)$ is denoted by $\chi^h$.  By \cite[Theorem 2]{CappellFree},
 $\Split(-\Split(h)\cdot h) = -\Split(h) + \Split(h) = 0$, so the first component of $\nncc$ is indeed splittable.

 In dimension $n=5$, since the submanifold $S^4$ of $\RP^5 \# \RP^5$ defining the connected sum is simply-connected, the bijection $\nncc$ exists
by the modification in \cite[Theorem 5, \S V.2]{CappellSplitting} of the $4$-manifold stable surgery in \cite[Theorems 4.1, 5.1]{CS}.  Alternatively, proceed analogous to the $n = 4$ case below.

In dimension $n=4$, we prove the existence of the bijection $\nncc$ by the following indirect method. For $h \in \cS(\RP^4\#\RP^4)$ one can, following Cappell, define the obstruction  $
\Split(h) \in \unil_5(\Z;\Z^-,\Z^-)$. Since the fundamental group $D_\infty=\Z\rtimes\Z_2$ is small, by Wall realization and plus construction
\cite[Theorems 11.3A, 11.1A]{FQ}, there exists a topological normal bordism $H$ from $h$ to another homotopy equivalence $h' = -\Split(h)\cdot h$ with surgery obstruction
$\sigma_*(H)=i(-\Split(h))$, where $i : \unil_5(\Z;\Z^-,\Z^-) \hookrightarrow L_5(\Z[\Z_2^- * \Z_2^-])$.    Furthermore $\Split(h') = 0$ as before.  It remains to show that $h'$ is splittable.  We accomplish this by periodicity.  (See \cite[Theorem 9.9]{Wall}).

Let $Z$ be a closed, simply-connected $4k$-manifold with signature $\sigma(Z)$.
$$
\Split(h' \times \id_Z) = \sigma(Z) \Split(h') =0.
$$
Taking $Z = P_2(\C)$ and using \cite[Proof of Theorem 1]{Weinberger}, one sees
that $h'$ is $\Z$-homology splittable, i.e. $h'$ is homotopic to a map whose restriction to the transverse inverse image $\Sigma^3$ of the
connecting $S^3$ in $\RP^4\#\RP^4$ is a $\Z$-homology equivalence. But then $h'$ is necessarily (topologically)
splittable, by the ``neck exchange'' trick of Jahren-Kwasik \cite[Proof of Theorem 2]{JK}--essentially due to the fact \cite[Cor. 9.3C]{FQ} that the
homology $3$-sphere $\Sigma$ bounds a contractible $4$-manifold $\Delta$.

\noindent(2) Let $h : M \to \RP^n\#\RP^n$ represent an element of the structure set.

By Lemma \ref{gen}, every element of $\haut(\RP^n\#\RP^n)$ can be represented by a map $\gamma$ which is the identity on the codimension one sphere $S$ defining the connected sum and which sends the complement $\RP^n\#\RP^n - S$ to itself (possibly interchanging the two components).  It is then clear that $\Split(\gamma \circ h) = \gamma_* \Split(h)$ where $\gamma_*$ is either the switch map or the identity.

 Wall realization  provides a normal bordism $H$ from $h$ to $h'\in \cS_{\Split}(\RP^n\#\RP^n)$ with surgery obstruction $\sigma_*(H)=-i(\Split(h)) \in L_{n+1}(\Z[\Z_2^\epsilon *
\Z_2^\epsilon])$. Then $(\gamma \times \id_{[0,1]}) \circ H$ is a normal bordism from $\gamma \circ h$ to $\gamma \circ h'$ and, since $\gamma$ is a homotopy equivalence preserving the orientation character, the surgery obstruction satisfies
$$
\sigma_*((\gamma \times \id_{[0,1]}) \circ H) = \gamma_*(\sigma_*(H)) = i(-\Split(\gamma \circ h)).
$$
It follows that
\begin{align*}
\gamma \circ (-\Split(h) \cdot h) &= \gamma \circ h' \\
& = -\Split(\gamma \circ h) \cdot (\gamma \circ h).
\end{align*}

In other words, the first component of the map $\nncc$ is also equivariant with respect to the action of $\haut(P^n \# P^n)$.

\noindent(3)  Observe that the connected sum operation
\[
\#: \cS(\RP^n)\times \cS(\RP^n) \to \cS_{\Split}(\RP^n\# \RP^n)\]
 on
topological structure sets is surjective.  For $n > 4$ this  uses the affirmative solution to the Poincar\'e conjecture for homotopy spheres of dimension $n-1 > 3$, and for $n=4$ this uses the  Jahren-Kwasik neck exchange trick, replacing a homotopy sphere $\Sigma^3$ by a genuine $S^3$.

To show injectivity of $\#$ for all $n>3$, suppose that $h_1 \# h_2$ is $s$-bordant to $h_1' \# h_2'$ via $H$,
where each $h_i, h_i' \in \cS(\RP^n)$. By taking the inverse image of the connecting $S^{n-1}$ under a homeomorphism, we may assume that $H$ is
a homotopy. Then we have a relative splitting problem of the $(n+1)$-dimensional homotopy equivalence $H$ along $S^{n-1}\times [0,1]$, which is
already split along its boundary $S^{n-1}\times\{0,1\}$. The relative form of the nilpotent normal cobordism construction in particular provides
an $s$-bordism $H'$ from $h_1\# h_2$ to $h_1'\# h_2'$, which splits as a boundary connected sum $H'=H'_1 \natural H'_2$ along $S^{n-1}\times
[0,1]$. Hence $h_i, h_i'$ represent the same element in $\cS(\RP^n)$ for all $i=1,2$.

\end{proof}

\begin{proof}[Proof of Theorem \ref{main}]  Lemma \ref{Lemma_Reduction} and Lemma \ref{gen} show that to compute $\bar J_n$, the set of homeomorphism classes of closed manifolds homotopy equivalent to $P^n \# P^n$, one needs to compute the action of $\gamma_1, \gamma_2, \gamma_3$ on $\cS(P^n \# P^n)$.  In terms of the bijection of Lemma \ref{Lemma_SplitAction},
$$
\cS(P^n \# P^n) \leftrightarrow \cS(P^n) \times \cS(P^n) \times \unil_{n+1}(\Z;\Z^\epsilon,\Z^\epsilon),
$$
$\gamma_1 \cdot (h_1, h_2, x) = (h_2, h_1, \sw x)$.  (More on $\sw x$ later in this paper.)

Note that post-composition with $\gamma_3$ acts via the identity on $\cS_{\Split}(P^n \# P^n)$ since such a structure can be represented by a connected sum $h_1 \# h_2 : M_1 \# M_2 \to P^n \# P^n$ and $\gamma_3 \circ h_1 \# h_2  = h_1 \# h_2  \circ  \gamma$ where $\gamma$ is the homeomorphism given by a Dehn twist in the domain.  Also, $\gamma_3$ acts trivially on $ \unil_{n+1}(\Z;\Z^\epsilon,\Z^\epsilon)$ since it induces the identity on the fundamental group.  Thus $\gamma_3$ acts trivially on $\cS(P^n \# P^n)$ by Lemma \ref{Lemma_SplitAction}, parts (1) and (2).

Now assume that $n$ is odd.  Post-composition by reflection through
$P^{n-1}$ acts on $\cS(P^n)$ by multiplication by $h \mapsto
\Phi^{-1}(-\Phi(h))$ according to Theorem \ref{P_classify}.   Lemma
\ref{Lemma_SplitAction} then shows that $\gamma_2 \cdot (h_1,h_2,x)
\mapsto (\Phi^{-1}(-\Phi(h_1)),\Phi^{-1}(-\Phi(h_2)),x)$.

Theorem \ref{main} follows.
\end{proof}


The proof shows a bit more.  It shows that the quotient function\break $\cS(\RP ^n \# \RP^n)/\langle \sw \rangle \to J_n$ is a bijection.  It allows us to be more precise about
the map $ J_n \to \bar J_n$ when $n \equiv 3 \pmod 4$.  Using the notation of Theorem \ref{P_classify} and Theorem \ref{main}, the
elements $(\{\Phi^{-1}(y),\Phi^{-1}(z) \}, \vartheta)\in J_n$ and $(\{\Phi^{-1}(-y),\Phi^{-1}(-z) \}, \vartheta)\in J_n$ map to the same element
of $\bar J_n$.  This is the only way a pair of elements in $J_n$ can have the same image in $\bar J_n$.

\section{Switch action on $\unil_2(\Z;\Z,\Z)$} \label{unil_2}

\begin{theorem}\label{Thm_UNil2}
The switch map $\sw$ operates as the identity on $\unil_2(\Z;\Z,\Z)$.
\end{theorem}

\begin{lemma}\label{Lemma_Induced}
Let $r:\unil_n(\Z;\Z,\Z)\to NL_n(\Z)$ be the Connolly-Ranicki isomorphism and $i:\unil_n(\Z;\Z,\Z)\to L_n(\Z[D_\infty])$ be the injection defined
by Cappell.  Let  $F:NL_n(\Z)\to L_n(\Z[D_\infty])$ be the map  $F =i \circ r^{-1}$.  Then
$F$ is given by the formula:
\[
F[(E,\chi)] = [(\Z[D_\infty] \otimes_{\Z[t]} E, a\otimes\chi)]
\]
where ``$a$'' is the left $\Z[D_\infty]$-module endomorphism of $\Z[D_\infty]$ defined by right multiplication by $a \in D_\infty$.
\end{lemma}

\begin{proof}  It is easy to verify
\begin{align*}
G : NL_n(\Z) & \to L_n(\Z[D_\infty])\\
[(E,\chi)] & \mapsto [(\Z[D_\infty] \otimes_{\Z[t]} E, a\otimes\chi)]
\end{align*}
is a well-defined homomorphism.  Here $(E,\chi)$ is a quadratic Poincar\'e complex over $\Z[t]$  in the sense of Ranicki \cite{RanickiExact}.

We need to show $G \circ r = i$.
 We first check first this on $\unil_2(\Z;\Z,\Z)$.
Let $x=\left[(P,Q,\rho_1 a,\rho_2 b,\mu_1 a,\mu_2 b)\right]\in\unil_2(\Z;\Z,\Z)$ (see \cite{CappellUnitary} for the definition of the UNil groups).
From Connolly-Ranicki,
$$r(x) =\left[ \prn{ P[t]\oplus Q[t] , \SmMatrix{ \rho_1 & 1 \\ -1 & \rho_2 t } , \SmMatrix{ \mu_1 \\ \mu_2 t } }\right].$$
Then
\begin{align*}
 G(r(x)) &= \left[\prn{ \Z[D_\infty]\otimes_{\Z[t]} (P\oplus Q) , \SmMatrix{ \rho_1 a & a \\ -a & \rho_2 b } , \SmMatrix{ \mu_1 a \\ \mu_2 b } } \right]\\
 & = \left[\prn{ \Z[D_\infty]\otimes_{\Z[t]} (P \oplus Q) , \SmMatrix{ \rho_1 a & 1 \\ -1 & \rho_2 b } , \SmMatrix{ \mu_1 a \\ \mu_2 b } }\right]\\
 & = i(x).
 \end{align*}
 where the second equality is seen by pulling back along the isomorphism $\SmMatrix{a & 0 \\
 0 & 1}$ and the last equality is
the definition of $i(x)$ (see \cite{CappellUnitary}).

To deal with the odd-dimensional case, in order to appeal to Connolly-Ranicki we use the fact that all of the maps $G$, $R$, and $i$ can be defined for any ring $R$, not just $\Z$ and are functorial with respect to maps of rings.  Therefore the standard Shaneson splitting
argument of Connolly-Ranicki \cite[Prop. 19]{CR} implies that $G\circ r=i$ on $\unil_3(\Z;\Z,\Z)$ also.
\end{proof}

\begin{remark}
The map $F$ has a geometric interpretation which led us to the formula of the previous lemma.  A homotopy equivalence $X^n \to P^n \# P^n$ can be considered both as a one-sided splitting problem by splitting along $P^{n-1} \# P^{n-1}$ and a two-sided splitting problem by splitting along $S^{n-1}$.  The passage from the one-sided splitting obstruction to the two-sided splitting obstruction coincides with the map $F$.  More precisely  $F$  is the composite
\[
NL_n(\Z) \xra{} L_n(\Z[t]) \xra{} L_n(\Z[t,t\inv]) \xra{\aqk\inv} LN_n(\Z \to \Z_2^- * \Z_2^-) \xra{\partial} L_n(\Z[\Z_2*\Z_2]).
\]
The first map is an inclusion, and the second is induced by an inclusion of rings (note we have the trivial involution $t \mapsto t$).

The group $LN_n(\Z \to \Z_2^- * \Z_2^-)$ is defined in \cite[Chapter 11]{Wall} as the obstruction group for one-sided splitting and sits in the
exact sequence
$$
\cdots \to L_{n+1}(\Z[\Z_2^- * \Z_2^-]) \to L_{n+2}(\Z[\Z] \to \Z[\Z_2 * \Z_2])  \to LN_n(\Z \to \Z_2^- * \Z_2^-) \to \cdots
$$
where the first map is a transfer map corresponding to the line bundle of the  ``two-fold cover'' $B\Z \to BD_\infty$. This $LN$-group is
identified algebraically in \cite[\S 7.6 pp. 691--695]{RanickiExact} as a relative term in the exact sequence. That is, each element of $LN_n$
is represented by a pair consisting of a quadratic Poincar\'e complex of dimension $n$ over $\Z[\Z_2^- * \Z_2^-]$ and a quadratic pair of
dimension $n+2$ over $\Z[\Z] \to \Z[\Z_2*\Z_2]$ with algebraic boundary the transfer pair of the complex.

The map $\partial$ is defined in \cite[\S 7.2 p. 565]{RanickiExact} as the boundary map in the transfer sequence; it forgets the quadratic pair.
Geometrically it corresponds to a transversal restriction of a one-sided splitting problem to a degree one normal map.

The map $\aqk$ is an isomorphism.  It was defined geometrically by Wall (see \cite[Theorem 12.9]{Wall}) as the
$\Z_2$-equivariant defect for handle exchanges in the middle dimension of a certain regular two-fold cover.
The map $\aqk$ was dubbed the antiquadratic kernel by Ranicki \cite[\S 7.6 pp. 698--699]{RanickiExact} and was given an algebraic definition.

Actually for us,  $\aqk^{-1}$ is the relevant map and is somewhat easier to define.  However, we will omit the definition here, and refer to \cite[Proof 7.6.3, p. 702]{RanickiExact} for the formula
\[
(\partial \circ \aqk\inv)\left[(E,\chi )\right] = \left[(\Z[D_\infty] \otimes_{\Z[t]} E, a \otimes \chi)\right].
\]

A careful reading of \cite[\S 7.6 pp. 737--745]{RanickiExact} shows that $\partial$ takes the one-sided splitting obstruction to the image of the two-sided splitting obstruction and that the composite map $F$ above satisfies $F \circ r = i$. \qed
\end{remark}

\begin{proof}[Proof of Theorem \ref{Thm_UNil2}]
Recall \cite[Theorem 4.6(2)]{CD} (see also \cite{CK}) that the map $L_2(\Z[t]) \to L_2(\F_2[t])$ is an isomorphism and the latter is detected by
restriction of the Arf invariant for the characteristic two field $\F_2(t)$ -- first map to $L_2(\F_2(t))$ and then by the Arf invariant to the
idempotent quotient $\F_2(t)/\{ f^2 -f \;|\; f \in \F_2(t) \}$. So it suffices to compute the switch map $\sw$ on the values of the inverse map
\[
\frac{t \F_2[t]}{\{f^2 - f \;|\; f \in \F_2[t]\} } \to NL_2(\Z)
\]
defined by
\[
[t p] \mapsto [P_{tp,1}] = \brk{\Z[t]^2, \SmMatrix{0 & 1\\ -1 & 0}, \SmMatrix{tp\\ 1}}.
 \]

Write $t = ba$, which generates the infinite cyclic subgroup $\Z$ of index two in $D_\infty$. Recall the definition of $F: NL_2(\Z) \to
L_2(\Z[D_\infty])$ in Lemma \ref{Lemma_Induced}.

Then $F[P_{tp,1}]$ is represented by the quadratic form:
$$\prn{ \Z[D_\infty]^2 , \SmMatrix{0 & a \\ -a & 0} , \SmMatrix{ p(t)b \\ a } }. $$
Note
\begin{align*}
\sw F[P_{tp,1}] &= \sw \brk{ \Z[D_\infty]^2 , \SmMatrix{0 & a \\ -a & 0} , \SmMatrix{ p(t)b \\ a } } \\
                  &= \sw \left\{ \SmMatrix{b & 0 \\ 0 & a}^* [ \Z[D_\infty]^2 , \SmMatrix{0 & a \\ -a & 0} , \SmMatrix{ p(t)b \\ a } ] \right\} \\
                                    &= \sw   [ \Z[D_\infty]^2 , \SmMatrix{b & 0 \\ 0 & a} \SmMatrix{0 & a \\ -a & 0} \SmMatrix{b & 0 \\ 0 & a} , \SmMatrix{b p(t)bb \\ aaa } ]  \\
                  &= \sw \brk{ \Z[D_\infty]^2 , \SmMatrix{ 0 & b \\ -b & 0} , \SmMatrix{  bp(t) \\ a} } \\
                  &= \brk{ \Z[D_\infty]^2 , \SmMatrix{ 0 & a \\ -a & 0} , \SmMatrix{ ap(t^{-1}) \\ b} } \\
                  &= \brk{ \Z[D_\infty]^2 , \SmMatrix{ 0 & a \\ -a & 0} , \SmMatrix{ p(t)a \\ ta} } \\
                  &= F[P_{p,t}].
                  \end{align*}
Therefore $\sw[P_{tp,1}]=[P_{p,t}]$ by injectivity of $F$. But both
$[P_{p,t}]$ and $[P_{tp,1}]$ have the same Arf invariant $[t p]$.
\end{proof}


\section{Switch action on $\unil_3(\Z;\Z,\Z)$} \label{unil_3}

We compute the switch action on $\unil_3(\Z;\Z,\Z)$ in a similar way to what we did for $\unil_2(\Z;\Z,\Z)$; we use \cite{CD} to find generators, use Lemma \ref{Lemma_Induced} to compute the switch map on the generators, and use the obstruction theory of \cite{CD} to express the result in terms of the original generators.


Connolly-Davis \cite{CD} give a $\Z$-module isomorphism
$$\SmMatrix{j_1 & j_2}:\frac{t\Z_4[t]}{\{2p(t^2)-2p(t) : p(t) \in t\Z_4[t] \}}\times t\Z_2[t]\to \unil_3(\Z;\Z,\Z).$$


\begin{theorem}\label{Theorem_Verschiebung}
Let $\pi:\frac{t\Z_4[t]}{\{2p(t^2)-2p(t)\}}\to t\Z_2[t]$ be the quotient map $\pi[tp] = [tp]$.
Then using the above coordinates $\SmMatrix{j_1 & j_2}$ for $\unil_3(\Z;\Z,\Z)$, the switch map is $\sw = \SmMatrix{1 & 0\\ \pi & 1}$.
\end{theorem}

\begin{corollary}
$x \in \unil_3(\Z;\Z,\Z)$ is fixed by the switch map if $x$ is divisible by two.
\end{corollary}

\begin{example}
Consider the Poincar\'e $(-1)$-quadratic complex $(C,\psi)$ of dimension one over $\Z[D_\infty]$ defined by the data:
\begin{itemize}
\item $C_1 = \Z[D_\infty] \oplus \Z[D_\infty] = C_0$ and $d := \SmMatrix{2 & 0\\ 0 & 2}$,%

\item $\psi_0 = \SmMatrix{2a & 1\\ 1 & b} : C^0 \to C_1$ and $\psi_1 = \SmMatrix{-2a & 0\\ -2 & -b}: C^0 \to
C_0$.
\end{itemize}
Observe it is the image of a certain $(-1)$-quadratic nilcomplex of dimension one over the triple $(\Z; \Z, \Z)$, see \cite[Definition
11.11]{Brookman}, under the natural group monomorphism $\unil_3(\Z;\Z,\Z) \to L_3(\Z[D_\infty])$ of cobordism classes. Moreover from the proof of Lemma
\ref{Lemma_Linking}, its cobordism class $[(C,\psi)]$ is also the image of an order $4$ element $[\cN_{t,1}] = j_1[t]$ of  $t\Z_4[t]/\{2p(t^2)-2p(t)\}$. Then
\[
\sw[\cN_{t,1}] = \sw(j_1[t]\oplus 0) = j_1[t]\oplus j_2[t] = [\cN_{1,t}].
\]
In particular $\sw[(C,\psi)] \neq \pm[(C,\psi)]$, hence $[(C,\psi)] \in L_3(\Z[D_\infty])$ is the image of a non-zero element of $\unil_3(\Z;\Z,\Z)$
which is not fixed by the switch map.
\end{example}


The actual computation of Connolly-Davis was of the group $\wt{\cL}(\Z[t],2)$, the Witt group of quadratic linking forms on $\Z[t]$-modules with exponent 2 which become Witt trivial under the map $t \to 0$.  Recall that these are defined to be triples $(M,b,q)$ where
\begin{itemize}
\item $M$ is a $\Z[t]$-module isomorphic to $\F_2[t]^k$ for some $k$;
\item $b:M\times M\to \Q[t]/\Z[t]$ is a nonsingular symmetric linking form;
\item $q:M\times M\to \Q[t]/2\Z[t]$ is a quadratic refinement of $b$ (so that in particular $[q(x)]=b(x,x)\in\Q[t]/\Z[t]$).
\end{itemize}

This group was shown to be isomorphic to $\unil_3(\Z;\Z,\Z)$ via the composite isomorphism $$\unil_3(\Z;\Z,\Z)\xrightarrow{\cong} NL_3(\Z)\xleftarrow{\cong} NL_0(\Z,\langle 2\rangle)\xleftarrow{\cong} \wt{\cL}(\Z[t],2).$$  Here the first isomorphism was from Connolly-Ranicki \cite{CR}, the second from a localization exact sequence, and the third by a devissage argument.

The maps $j_1$ and $j_2$ were defined in terms of certain basic linking forms, so it is sufficient for us to compute the effect of the switch map on those:

\begin{definition}[{\cite[Definition 1.6]{CD}}]
For $p,g\in\Z[t]$ where either $p(0) = 0$ or $q(0) = 0$, define $[\cN_{p,g}]\in\wt{\cL}(\Z[t],2)$ by the $(+1)$-quadratic linking form
$$\cN_{p,g}=\left( \Z_2[t]^2, \begin{pmatrix}p/2 & 1/2 \\ 1/2 & 0 \end{pmatrix}, \begin{pmatrix} p/2 \\ g \end{pmatrix} \right).$$
\end{definition}

Then $j_1[tp] = \left[\cN_{tp,1}\right]$ and $j_2[tp] = \left[\cN_{1,tp}\right] - \left[\cN_{t,p}\right]$.

\begin{lemma}\label{Lemma_Complexes}
Let $s : \wt{\cL}(\Z[t],2) \to NL_3(\Z)$ be the isomorphism mentioned above.
\begin{enumerate}
\item $s[\cN_{tp,g}]$ is represented by the $(-1)$-quadratic 1-dimensional Poincar\'e complex in $NL_3(\Z)$:
$$\xymatrix@+2mm{ C^0=\Z[t]^{2} \ar[r]^2 \ar[d]_{\SmMatrix{p(t)t & 1 \\ 1 & 2g(t)}} \ar[dr]|{\SmMatrix{-p(t)t & -1 \\ -1 & -2g(t)}} & \Z[t]^{2}=C^1\ar[d]^0 \\
             C_1=\Z[t]^{2} \ar[r]_2 & \Z[t]^{2}=C_0 } $$
\item $F(s[\cN_{tp,g}])$ is represented by the $(-1)$-quadratic 1-dimensional Poincar\'e complex in $L_3(\Z[D_\infty])$:
$$\xymatrix@+2mm{ C^0=\Z[D_\infty]^{2} \ar[r]^2 \ar[d]_{\SmMatrix{p(t)b & a \\ a & 2g(t)a}} \ar[dr]|{\SmMatrix{-p(t)b & -a \\ -a & -2g(t)a}} & \Z[D_\infty]^{2}=C^1\ar[d]^0 \\
             C_1=\Z[D_\infty]^{2} \ar[r]_2 & \Z[D_\infty]^{2}=C_0 } $$
\end{enumerate}

\end{lemma}
\begin{proof}
The isomorphism $\wt\cL(\Z[t],2)\cong NL_0(\Z,\langle 2\rangle)$ is the obvious map $[(M,b,q)]\to [(M,b,q)]$.

Ranicki \cite[Prop. 3.4.1]{RanickiExact} showed that $NL_0(\Z,\langle 2\rangle)$ can equivalently be described as cobordism classes of $\langle
2\rangle$-acyclic 1-dimensional $(-1)$-quadratic Poincar\'e complexes over $\Z[t]$.  The boundary map $NL_0(\Z,\langle 2\rangle)\to NL_3(\Z)$ in
the localization exact sequence is then given by the map $[(C,\psi)]\to [(C,\psi)]$. The content of the first part of the lemma is thus to find the
1-dimensional Poincar\'e complexes corresponding to the linking forms $\cN_{tp,g}$.
It follows immediately by applying the correspondence from \cite{RanickiExact}, recalled below in Definition \ref{Lemma_Resolutions}, that $\cN_{tp,g}$ corresponds to the quadratic Poincar\'e complex claimed.

The second part is then an immediate consequence of Lemma \ref{Lemma_Induced}.  (Recall that $p(t)ta=p(t)baa=p(t)b$.)
\end{proof}

\begin{definition}[{\cite[Proof 3.4.1]{RanickiExact}}]\label{Lemma_Resolutions}
Let $R$ be a ring with involution, and let $S$ be a central multiplicative subset of $R$ which is invariant under the involution. Let $(C,\psi)$
be a $S$-acyclic 1-dimensional $(-1)$-quadratic Poincar\'e complex over $R$, so that $(C,\psi)$ consists of:
\begin{itemize}
\item a monomorphism $d: C_1 \to C_0$ of finitely generated free $R$-modules such that $\mathrm{id}_{S\inv R} \otimes_R d$ is
an isomorphism,
\item a morphism $\psi_0: C^0 \to C_1$ such that $\psi_0^*: C^1 \to C_0$ induces an isomorphism $\mathrm{Cok}(d^*) \to \mathrm{Cok}(d)$, and
\item a morphism $\psi_1: C^0 \to C_0$ such that $\psi_1 + \psi_1^* = - d \circ \psi_0$.
\end{itemize}

{\em The nonsingular $(+1)$-quadratic linking form $(M,b,q)$ over $(R,S)$ associated to the resolution $(C,\psi)$} consists of:
\begin{itemize}
\item the $S$-torsion $R$-module $M := \mathrm{Cok}(d^*)$ of homological
dimension one,
\item the sesquilinear map $b: M \times M \to S\inv R/R$, with adjoint
an isomorphism, well defined for all representatives $x, y \in C^1$ by $b(x,y) := \frac{1}{s} \Inn{y}{\psi_0(z)}$ where $z \in C^0$ is uniquely
determined by the formula $d^*(z) = s x$ for any given $s \in S$ such that $s [x] = 0 \in M$, and
\item the quadratic map $q: M \to S\inv R/\{r + \ol{r} \;|\; r \in
R\}$, well defined for all $x$ denoted as above by $q(x) := \frac{1}{s} \Inn{z}{\psi_1(z)}$.
\end{itemize}
\end{definition}

We are now in a position to compute the switch involution on the generators of $\wt\cL(\Z[t],2)$.

\begin{lemma}\label{Lemma_Linking}
For all polynomials $p,g \in \Z[t]$, the involution $\sw$ is given by:
\[
\sw[\cN_{tp,g}] = [\cN_{p,tg}].
\]
\end{lemma}

\begin{proof}
Denote the Poincar\'e complex over $\Z[D_\infty]$ of Lemma
\ref{Lemma_Complexes}(2) by $\SmBMatrix{p(t) b & a\\ a & 2 g(t) a}$.
Similar to the proof of Theorem \ref{Thm_UNil2}, note that
\begin{align*}
\sw\SmBMatrix{p(t) b & a\\ a & 2 g(t) a} &= \sw \left\{\SmMatrix{b & 0 \\ 0 & a}^*\SmBMatrix{p(t) b & a\\ a & 2 g(t) a}\right\}\\
 &= \sw \left\{ \SmMatrix{b & 0 \\ 0 & a}\SmBMatrix{p(t) b & a\\ a & 2 g(t) a}\SmMatrix{b & 0 \\ 0 & a}\right\}\\
&= \sw \SmBMatrix{bp(t) & b \\ b & 2ag(t)} \\
&= \SmBMatrix{ap(t^{-1}) & a \\ a & 2bg(t^{-1})} \\
&= \SmBMatrix{p(t)a & a \\ a & 2g(t)baa}
 \end{align*}
This is precisely $F(s[\cN_{p,tg}])$.  Hence $\sw F(s[\cN_{tp,g}])=F(s[\cN_{p,tg}])$, so by injectivity of $F$ and $s$, $\sw[\cN_{tp,g}]=[\cN_{p,tg}]$.
\end{proof}


\begin{proof}[Proof of Theorem \ref{Theorem_Verschiebung}]
Let $p,p' \in \Z[t]$ be polynomials. Recall from \cite[p. 1072]{CD} that $j_1[tp]:=[\cN_{tp,1}]$, and that $j_2[tp']:= [\cN_{1,tp'}] - [\cN_{t,p'}]$.
Note by Lemma
\ref{Lemma_Linking} that:
\begin{align*}
\sw(j_1[tp]) &= \sw([\cN_{tp,1}])\\
 &= [\cN_{p,t}]\\
&= ([\cN_{p,t}] - [\cN_{1,tp}] + [\cN_{t,p}]) \oplus j_2[tp]\\
\sw(j_2[tp']) &= \sw([\cN_{1,tp'}] - [\cN_{t,p'}])\\
 &= [\cN_{t,p'}]- [\cN_{1,tp'}]\\
&= 0\oplus j_2[tp'].
 \end{align*}
It remains to show $[\cN_{t,p}] + [\cN_{p,t}] = [\cN_{1,tp}] +
[\cN_{tp,1}]$.

Recall
\begin{align*}
\cN_{t,p}\oplus \cN_{p,t} &=
 \prn{\F_2[t]^4, \SmMatrix{t/2 & 1/2 & 0 & 0\\ 1/2 & 0 & 0 & 0\\ 0 & 0 & p/2 & 1/2\\ 0 & 0 & 1/2 & 0}, \SmMatrix{t/2\\ p\\ p/2\\ t}}\\
\cN_{1,tp}\oplus \cN_{tp,1} &=
 \prn{\F_2[t]^4, \SmMatrix{1/2 & 1/2 & 0 & 0\\ 1/2 & 0 & 0 & 0\\ 0 & 0 & tp/2 & 1/2\\ 0 & 0 & 1/2 & 0}, \SmMatrix{1/2\\ tp\\ tp/2\\
 1}}.
\end{align*}

For a polynomial $p \in \F_2[t]$, define polynomials $p_\even, p_\odd \in \F_2[t]$ by the equation $p = p^2_\even + t p^2_\odd$.
Label the basis of $\F_2[t]^8$ as $\{e_1, e_2, \dots , e_8\}$ and  let
\begin{align*}
v_0 &= p_\even \cdot e_4 + e_6 + t p_\odd \cdot e_8\\
v_1 &= e_2 + p_\odd\cdot e_4 + p_\even\cdot e_8.
\end{align*}
Then $S := \Span{v_0,v_1}$ is a sublagrangian of  the $(+1)$-quadratic linking form $(\cN_{t,p}\oplus \cN_{p,t}) - (\cN_{1,tp}\oplus
\cN_{tp,1})$ of exponent $2$ over $(\Z,\langle 2\rangle)$.
Note
\[
S^\perp  = \Span{ p_\odd\cdot e_1 + e_3 + p_\even \cdot e_5,
e_4,  p_\even\cdot e_1  + t p_\odd \cdot e_5+ e_7 , e_8, v_0, v_1}.
\]
Then $(\cN_{t,p}\oplus \cN_{p,t}) - (\cN_{1,tp}\oplus \cN_{tp,1})$
is Witt equivalent to the sublagrangian construction:
\[
(S^\perp/S, \ol{b}, \ol{q}) = \prn{\F_2[t]^4, \SmMatrix{0 & 1/2 & 0
& 0\\ 1/2 & 0 & 0 & 0\\ 0 & 0 & 0 & 1/2\\ 0 & 0 & 1/2 & 0}, \SmMatrix{p\\ t\\ tp\\
1}}.
\]
Lemma 4.3(2) and (6) of \cite{CD} show that this admits a lagrangian.  Alternatively, note that the above is an even quadratic linking form so
represents an element of   $\wt{\cL}^\even(\Z[t],2)$, the Witt group of even quadratic linking forms $(M,b,q)$ on $\Z[t]$-modules with exponent
2. Here $(M,b,q)$ is  {\em even} if $b(x,x) \in \Z[t]$ for all $x \in M$.  By an easy change of coordinates, one sees that
$\wt{\cL}^\even(\Z[t],2) \cong NL_0(\F_2)$ which is isomorphic to ${t \F_2[t]}/{\{f^2 - f \;|\; f \in \F_2[t]\} }$ via the Arf invariant of the
function field $\F_2(t)$.  On the other hand, the Arf invariant of the above form is trivial, hence the form is Witt trivial.

\end{proof}

The above proof is complete, but the last part of the proof was unmotivated.  The motivation is the obstruction theory of Connolly-Davis.  There
is a short exact sequence of $\Z$-modules, which first appeared in \cite[Theorem 25]{CR},
$$
0 \to \wt{\cL}^\even(\Z[t],2) \to \wt{\cL}(\Z[t],2) \xrightarrow{B} t \Z_2[t] \oplus t \Z_2[t] \to 0.
$$
See \cite[Definition 6.1]{CD} for the definition of $B$ and \cite[Lemma 5.7(1), Example 6.2]{CD} for exactness.

 Consider $B$ as a ``characteristic number'' -- a certain combination of quadratic values of the Wu classes \cite[Equation 5.(4)]{CD} $v_0, v_1$ -- which is the obstruction \cite[\S 6 Proof 1.7 Line 3]{CD} for a quadratic
linking form to Witt equivalent to even form.  If the quadratic linking form is even,  the Arf invariant over the characteristic two field $\F_2(t)$ can then be applied.

\begin{corollary}
The switch map on the invariant $B$ is: $B \circ \sw = \SmMatrix{1 & 0\\ 1 & 1}\circ B$.
\end{corollary}

\begin{proof}
Recall from \cite[Equations (13),(14)]{CD} that:
\[\begin{array}{lllllll}
B_1 \circ j_1 & = & \pi & \qquad & B_1 \circ j_2 & = & 0\\
B_2 \circ j_1 & = & 0 & \qquad & B_2 \circ j_2 & = & \id.
\end{array}\]
Then for any polynomials $p, p' \in \Z[t]$ note
\[
(B\circ \sw)(j_1[tp] + j_2[tp']) \; =\; \begin{pmatrix}B_1\\
B_2\end{pmatrix}
\prn{j_1[tp] + j_2[tp] + j_2[tp']} \; =\; \begin{pmatrix}[tp]\\
[tp] + [tp']\end{pmatrix}.
\]
\end{proof}


\bibliographystyle{alpha}

\end{document}